\documentclass{amsart}
\usepackage{amsmath,amssymb,amsthm,mathtools,amsfonts}
\usepackage{booktabs,longtable,array}
\usepackage{multirow}
\usepackage{graphicx}
\usepackage{capt-of}
\usepackage{tikz}
\usetikzlibrary{patterns}
\usetikzlibrary{decorations.markings}
\usetikzlibrary{calc}
\newcommand{\F}{\mathbb F}
\newcommand{\PG}{\operatorname{PG}}
\newcommand{\Span}[1]{\langle #1\rangle}
\newcommand{\cS}{\mathcal S}
\newcommand{\supp}{\operatorname{supp}}
\newcommand{\wt}{\operatorname{wt}}
\newcommand{\diam}{\operatorname{diam}}

\newtheorem{theorem}{Theorem}
\newtheorem{prop}{Proposition}
\newtheorem{lem}{Lemma}

\newtheorem{rem}{Remark}

\numberwithin{equation}{section}

\title[Distances in the collinearity graphs of simplex codes]{Distances in the collinearity graphs of geometries related to simplex codes}
\author{Adam Tyc}

\keywords{simplex code, point-line geometry, collinearity graph, graph diameter, graph distance}
\subjclass[2020]{51E20, 51E22, 05C12}

\address{Faculty of Mathematics and Computer Science, 
	University of Warmia and Mazury, S{\l}oneczna 54, 10-710 Olsztyn, Poland}
\email{adam.tyc@uwm.edu.pl}

\begin{document}
	\begin{abstract}
		Let $\cS(k,q)$ be the point-line geometry whose points and lines correspond to one-dimensional and two-dimensional subcodes of $q$-ary simplex codes of dimension $k$, respectively. We determine the complete distance relation and the diameter of its
		collinearity graph. For two points whose representing vectors have different sets of zero coordinates there is a common neighbour. If the representing vectors have the same zero set, then a common neighbour exists precisely when $(q-1)q^{k-2}$ of their non-zero coordinate ratios can be selected so that their product is one. To prove this, we use a generalization of Hall's result for finite abelian groups. Applying the Erd\H{o}s--Ginzburg--Ziv theorem, we establish that the above product condition holds for $k\ge3$, which implies that the diameter is two. If $k=2$, then the diameter is one when $q=2$, two when $q=3$ and three when $q\ge4$.
	\end{abstract}
	\maketitle
	
	\section{Introduction}
	
	Graphs and geometries related to various classes of linear codes are studied in \cite{CGK21,CG24,KP16,KPP18,KP20,KP24,KPT26,PPZ25,Siem26}.
	In every such class, linear codes are of the same length and dimension.
	Two codes are adjacent vertices of the graph if their intersection has codimension one in each of these codes. 
	The graphs of non-degenerate linear codes and projective codes are considered in 
	\cite{CGK21,CG24,KP16,KPP18}.
	
	Let $q$ be a prime power and let $k\ge2$. 
	Then $n=(q^k-1)/(q-1)$ is the number of points in $\PG(k-1,q)$. 
	The \emph{$q$-ary simplex codes of dimension $k$} can be characterized as the
	subspaces of $\F_q^n$ maximal with respect to the property that every non-zero
	codeword has Hamming weight $q^{k-1}$ \cite{HP03}. 
	
	Consider the graph $\Gamma^s(k,q)$ whose vertices are the $q$-ary simplex codes of
	dimension $k$ (as above two such codes are adjacent if their intersection has dimension $k-1$). 
	In \cite{KPP18} this graph is presented as a special case
	of the graph of projective codes. 
	A detailed study of $\Gamma^s(2,4)$ is given in 
	\cite{KP20}. Maximal cliques of $\Gamma^s(k,q)$ are investigated in \cite{KP24}.
	The diameter of $\Gamma^s(2,q)$ is determined in \cite{Siem26}.
	
	The related point-line geometry is considered in \cite[Section~6]{KP24}.
	Let $\cS(k,q)$ be the geometry whose points and lines are the
	one-dimensional and two-dimensional subspaces of $\F_q^n$ contained in simplex
	codes of dimension $k$, respectively. The geometry $\cS(k,2)$ also occurs in the more general
	family of point-line geometries related to binary equidistant codes studied in
	\cite{PPZ25}. We denote the \emph{collinearity graph} of $\cS(k,q)$ by
	$\Delta(k,q)$. Its vertices are the points of $\cS(k,q)$ and two distinct
	vertices are adjacent if they are joined by a line in $\cS(k,q)$. Maximal
	cliques in this graph are investigated in \cite{KPT26}; 
	note that some
	maximal cliques of $\Delta(2,q)$ for $q\ge5$ are normal rational curves.
	
	In the present paper, we determine the complete distance relation in $\Delta(k,q)$.
	
	The vertices of $\Delta(k,q)$ are projective points represented by vectors of Hamming weight $q^{k-1}$. The collinearity relation of $\Delta(k,q)$ is more restrictive than this weight condition. 
	Consider representatives of two such points as the rows of a $(2\times n)$-matrix. 
	The collinearity of the corresponding points can be characterized by the numbers of zero columns and mutually proportional non-zero columns (Proposition \ref{prop:subcodes}). 
	Using this characterization, we express adjacency of $\Delta(k,q)$ in terms of zero coordinates and ratios of non-zero coordinates.
	
	The distance relation has a different form according to whether two points have the same set of zero coordinates.
	If these sets are different, any two non-adjacent points have a common neighbour. The proof of this fact is based on completing two lists of ratios such that every element of $\F_q^\times$ occurs $q^{k-2}$ times in each of them
	(a weak version of these arguments is used in \cite{Siem26} for the graph of two-dimensional simplex codes). 
	To obtain such a completion we generalize Hall's result \cite{Hall52} concerning permutations of finite abelian groups.
	
	When the sets of zero coordinates coincide, the existence of a common neighbour is determined by a product condition on the non-zero coordinate ratios. We show that this condition is both necessary and sufficient.
	For $k\ge3$, we use the well-known Erd\H{o}s--Ginzburg--Ziv theorem \cite{EGZ61} to show that this condition is always satisfied. This guarantees that the diameter of $\Delta(k,q)$ is two.
	If $k=2$, then the product condition is not necessarily satisfied;  
	when this happens the corresponding pair is at distance three.
	In this case, the diameter is three for $q\ge4$, while for $q=2$ and $q=3$ it is one and two, respectively.

	\section{Preliminaries}\label{sec:preliminaries}
	\subsection{Simplex codes and their subcodes}
	
	Let $q$ be a prime power and let $\F_q$ be the field with $q$ elements. A \emph{$q$-ary linear code of length $n$} is a vector subspace of $\F_q^n$. If its dimension is $k$, then it is called a \emph{linear $[n,k]_q$ code}. A matrix whose rows form a basis of the code is a \emph{generator matrix}.
	
	For $m\ge1$, let $\PG(m-1,q)$ denote the projective space of one-dimensional subspaces of $\F_q^m$. A non-zero vector $v\in\F_q^m$ represents the projective point $\langle v\rangle$; two non-zero vectors represent the same point if and only if they are proportional.
	
	For an integer $s\ge0$, write
	$$
	[s]_q=\frac{q^s-1}{q-1}.
	$$
	An $s$-dimensional vector space over $\F_q$ has $[s]_q$ one-dimensional subspaces, so $\PG(s-1,q)$ has $[s]_q$ points when $s\ge1$.
	
	Fix $k\ge2$. We set
	$$
	n=[k]_q,
	$$
	and let $V=\F_q^n$. A linear $[n,k]_q$ code is a \emph{simplex code} if the columns of a generator matrix are non-zero and mutually non-proportional. Since there are $n=[k]_q$ columns, they represent all points of $\PG(k-1,q)$. Changing the row basis acts on the columns by an invertible linear transformation, so this definition is independent of the generator matrix. A \emph{subcode} is a subspace of a code.
	
	We will use the following characterization of subcodes of simplex codes \cite[Proposition~1]{KP24}.
	
	\begin{prop}\label{prop:subcodes}
		Let $X\subset V$ be an $m$-dimensional subspace, where $1\le m\le k$, and let $M$ be a generator matrix of $X$. Then $X$ is contained in a simplex code of dimension $k$ if and only if $M$ has precisely $[k-m]_q$ zero columns and each point of $\PG(m-1,q)$ is represented by precisely $q^{k-m}$ non-zero columns of $M$.
	\end{prop}
	
	\subsection{Simplex points and collinearity}
	A point of $\PG(n-1,q)$ is a one-dimensional subspace
	$\langle x\rangle$ of $V$, where $x\in V\setminus\{0\}$.
	For $x=(x_1,\ldots,x_n)$, its \emph{support} is
	$\supp(x)=\{i:x_i\ne0\}$ and its \emph{Hamming weight} is
	$\wt(x)=|\supp(x)|$.
	
	We denote by $\cS(k,q)$ the geometry whose points are the points of
	$\PG(n-1,q)$ contained in simplex codes of dimension $k$.
	Its lines are the sets
	$$\{\langle z\rangle : z\in X\setminus\{0\}\},$$
	where $X$ is a two-dimensional subspace contained in a simplex code of
	dimension $k$.
	The points of $\cS(k,q)$ are called \emph{simplex points}.
	
	By Proposition~\ref{prop:subcodes} with $m=1$, a point $P=\Span{x}$ is a simplex point if and only if
	$$
	\wt(x)=n-[k-1]_q=q^{k-1}.
	$$
	Its \emph{zero set} is
	$$
	Z(P)=\{i\in\{1,\ldots,n\}:x_i=0\}.
	$$
	It is independent of the non-zero representative $x$ and has size $h=[k-1]_q$. 
	For two-dimensional subcodes, set
	$$
	c=[k-2]_q,\qquad \mu=q^{k-2}.
	$$
	By Proposition~\ref{prop:subcodes}, a generator matrix of such a subcode has $c$ zero columns and each point of $\PG(1,q)$ is represented by precisely $\mu$ non-zero columns.
	We have
	\begin{equation}\label{eq:parameters}
		h=c+\mu,\qquad n=c+(q+1)\mu,\qquad c=\frac{\mu-1}{q-1}<\mu.
	\end{equation}
	
	The \emph{collinearity graph} $\Delta(k,q)$ has the simplex points as its vertices. Two distinct vertices $P,Q$ are adjacent if the projective line joining them is a line of $\cS(k,q)$; in this case we write $P\sim Q$. We denote the graph distance in $\Delta(k,q)$ by $d_\Delta$.
	
	\begin{lem}\label{lem:ratio}
		Let $P=\Span{x}$ and $Q=\Span{y}$ be distinct simplex points. Define $D=\supp(x)\cap\supp(y)$. 
		Then $P\sim Q$ if and only if
		\begin{equation}\label{eq:adjacency}
			|Z(P)\cap Z(Q)|=c
			\quad\text{and}\quad
			\bigl|\{i\in D:y_i/x_i=\alpha\}\bigr|=\mu
			\quad\text{for every }\alpha\in\F_q^\times.
		\end{equation}
	\end{lem}
	\begin{proof}
		The vectors $x,y$ are independent, so they generate a two-dimensional subspace. The zero columns of its generator matrix with rows $x,y$ are precisely those indexed by $Z(P)\cap Z(Q)$. If there are $c$ of them, each of $Z(P)\setminus Z(Q)$ and $Z(Q)\setminus Z(P)$ has size $h-c=\mu$. The corresponding columns represent the points $[0:1]$ and $[1:0]$ of $\PG(1,q)$. The columns indexed by $D$ represent the remaining points $[1:\alpha]$, where $\alpha=y_i/x_i\ne0$. Proposition~\ref{prop:subcodes} with $m=2$ now gives the equivalence.
	\end{proof}
	
	Rescaling $x$ or $y$ multiplies all the non-zero ratios by a common element of $\F_q^\times$, which permutes that group. Thus condition~\eqref{eq:adjacency} is independent of the representatives. Since $h>c$, distinct points with the same zero set are never collinear. Every simplex point has a neighbour: it belongs to a simplex code of dimension $k\ge2$ and can be joined to any other projective point of that code.
	
	\section{Results}\label{sec:main}
	
	For distinct points $P=\Span{x}$ and $Q=\Span{y}$ of $\Delta(k,q)$ with the same zero set $Z$, we use the notation
	$$
	D=\{1,\ldots,n\}\setminus Z
	$$
	and for $i\in D$
	$$
	\lambda_i=\frac{y_i}{x_i},
	\qquad
	\Lambda=\prod_{i\in D}\lambda_i.
	$$
	There are $|D|=q\mu$ such ratios. The following theorem gives the complete
	distance characterization in $\Delta(k,q)$.
	
	\begin{theorem}\label{thm:distance}
		Let $P,Q$ be distinct simplex points of $\cS(k,q)$.
		\begin{enumerate}
			\item We have $d_\Delta(P,Q)=1$ if and only if \eqref{eq:adjacency} holds.
			\item If $P\not\sim Q$ and $Z(P)\ne Z(Q)$, then $d_\Delta(P,Q)=2$.
			\item If $Z(P)=Z(Q)$, then $d_\Delta(P,Q)=2$ if and only if there is a set of coordinate indices $F\subset D$ such that
			\begin{equation}\label{eq:product-condition}
				|F|=(q-1)\mu\qquad\text{and}\qquad\prod_{i\in F}\lambda_i=1.
			\end{equation}
			Equivalently, there is $W\subset D$ with $|W|=\mu$ and $\prod_{i\in W}\lambda_i=\Lambda$.
			If these equivalent conditions fail, then $d_\Delta(P,Q)=3$.
		\end{enumerate}
		In particular, $\Delta(k,q)$ is connected.
	\end{theorem}
	
	The condition in \eqref{eq:product-condition} is stated in terms of coordinate indices, so the same ratio value may occur more than once in the product if it appears at different indices. We also need to check that the condition does not depend on the chosen representatives of $P$ and $Q$. Replacing $y$ by $ty$ multiplies every ratio $\lambda_i$ by $t$ and hence multiplies the product over $F$ by $t^{|F|}=t^{(q-1)\mu}=1$. Replacing $x$ by a non-zero scalar multiple gives the same conclusion. Thus the product condition depends only on the projective points $P$ and $Q$.
	
	Using the Erd\H{o}s--Ginzburg--Ziv theorem 
	we establish that the condition \eqref{eq:product-condition} is always satisfied for $k\ge3$ (Section~\ref{sec:proofs}).
	For $k=2$, the condition reduces to
	$$
	\Lambda\in\{\lambda_i:i\in D\}
	$$
	(Section~\ref{sec:same}); 
	if this fails, then the distance is three.
	
	The following theorem determines the diameter of the collinearity graph.
	\begin{theorem}\label{thm:diameter}
		For every $k\ge2$ and every prime power $q$,
		$$
		\diam\Delta(k,q)=
		\begin{cases}
			1,&k=2,\ q=2,\\
			2,&k\ge3\text{ or }(k,q)=(2,3),\\
			3,&k=2,\ q\ge4.
		\end{cases}
		$$
	\end{theorem}
	
	\begin{rem}{\rm
			The case $(k,q)=(2,2)$ is trivial (there is only one simplex code in this case). 
			If $(k,q)=(3,2)$ or $(k,q)=(2,3)$, then $\cS(k,q)$ is a polar space of rank $3$ or $2$, respectively \cite{KPP18,KP20}. 
	}\end{rem}
	
	\section{Points with distinct zero sets}\label{sec:different}
	\subsection{A generalization of Hall's theorem}
	We first recall the result of Hall \cite{Hall52}. 
	
	\begin{theorem}[Hall \cite{Hall52}]\label{thm:Hall}
		Let $G$ be a finite (multiplicative) abelian group of order $N$ and let $b_1,\ldots,b_N\in G$ ($b_i$ need not be distinct). There exists a permutation $a_1,\ldots,a_N$ of the elements of $G$ such that $a_1b_1,\ldots,a_Nb_N$ is also a permutation of the elements of $G$ if and only if $\prod_{i=1}^N b_i=1$.
	\end{theorem}
	
	The above statement will be generalized as follows. 
	\begin{prop}\label{prop:Hall}
		Let $G$ be a finite (multiplicative) abelian group of order $r$ and let $m\ge1$ be an integer. Suppose that $\lambda_1,\ldots,\lambda_L\in G$ where $1\le L\le rm$.
		\begin{enumerate}
			\item If $L<rm$, there are $a_1,\ldots,a_L\in G$ such that every element of $G$ occurs at most $m$ times in each of the sequences
			\[
			a_1,\ldots,a_L
			\quad\text{and}\quad
			a_1\lambda_1^{-1},\ldots,a_L\lambda_L^{-1}.
			\]
			
			\item If $L=rm$, there are $a_1,\ldots,a_L\in G$ such that every element of $G$ occurs precisely $m$ times in each of the sequences
			\[
			a_1,\ldots,a_L
			\quad\text{and}\quad
			a_1\lambda_1^{-1},\ldots,a_L\lambda_L^{-1}
			\]
			if and only if
			\[
			\prod_{i=1}^L\lambda_i=1.
			\]
		\end{enumerate}
	\end{prop}
	
	\begin{proof}
		Let $H$ be a finite abelian group of order $m$ and consider $\widetilde G=G\times H$. 
		For $1\le i\le L$, let $b_i=(\lambda_i^{-1},1_H)$.
		
		Suppose first that $L<rm$. We extend the sequence $b_1,\ldots,b_L$ to a sequence
		$b_1,\ldots,b_{rm}$ by setting for $L<i<rm$
		$$b_i=1_{\widetilde G}$$
		and
		$$b_{rm}=\left(\prod_{i=1}^{L}b_i\right)^{-1}.$$
		Then
		$$\prod_{i=1}^{rm}b_i=1_{\widetilde G}.$$
		
		Apply Theorem~\ref{thm:Hall} to $\widetilde G$. It gives a permutation $A_1,\ldots,A_{rm}$ of the elements of $\widetilde G$ such that $A_1b_1,\ldots,A_{rm}b_{rm}$ is also a permutation of the elements of $\widetilde G$. Let $a_i$ be the first coordinate of $A_i$. Under the projection $\widetilde G\to G$, every element of $G$ has precisely $m$ preimages. Hence, among the first $L$ terms, every element of $G$ occurs at most $m$ times in each of the sequences
		$$a_1,\ldots,a_L
		\quad\text{and}\quad
		a_1\lambda_1^{-1},\ldots,a_L\lambda_L^{-1}.$$
		This proves (1).
		
		Suppose now that $L=rm$. If $\prod_{i=1}^L\lambda_i=1$, then
		$$\prod_{i=1}^L b_i=1_{\widetilde G}.$$
		Applying Theorem~\ref{thm:Hall} and projecting to $G$, we obtain $a_1,\ldots,a_L$ such that every element of $G$ occurs precisely $m$ times in each of the two sequences in (2).
		
		Conversely, suppose that such $a_1,\ldots,a_L$ exist. Since every element of $G$ occurs precisely $m$ times in both sequences, their products are equal. Hence
		$$\prod_{i=1}^L a_i
		=
		\prod_{i=1}^L a_i\lambda_i^{-1}.$$
		Cancelling $\prod_{i=1}^L a_i$ gives
		$$\prod_{i=1}^L\lambda_i=1.$$
	\end{proof}
	
	\subsection{Construction of a common neighbour}
	For $k=2$, the following result is proved in \cite[Lemma~2.4]{Siem26}; here it is extended to arbitrary $k$.
	
	\begin{lem}\label{lem:bridge}
		If two simplex points have different zero sets, then they have a common neighbour in $\Delta(k,q)$.
	\end{lem}
	\begin{proof}
		Let $P=\Span{x}$ and $Q=\Span{y}$. Partition the coordinate indices into
		$$
		\begin{aligned}
			A&=Z(P)\cap Z(Q),&
			B&=Z(P)\setminus Z(Q),\\
			C&=Z(Q)\setminus Z(P),&
			D&=\supp(x)\cap\supp(y).
		\end{aligned}
		$$
		Write $a=|A|$. Since the two zero sets have size $h$ and are different, $a<h$. Moreover,
		\begin{equation}\label{eq:blocks}
			|B|=|C|=h-a,\qquad |D|=n-2h+a=(q-1)\mu-c+a,
		\end{equation}
		where the last equality follows from \eqref{eq:parameters}.
		
		We first choose a set $E$ that will be the zero set of a common neighbour. It must have size $h$ and meet each of $Z(P),Z(Q)$ in precisely $c$ indices.
		Let $b=\max\{0,c-a\}$ and choose the indices of $E$ according to the following table:
		$$
		\begin{array}{c@{\qquad}c@{\qquad}c}
			\toprule
			\text{block}&\text{size}&\text{number of indices in }E\\
			\midrule
			A&a&c-b\\
			B&h-a&b\\
			C&h-a&b\\
			D&(q-1)\mu-c+a&\mu-b\\
			\bottomrule
		\end{array}
		$$
		All choices are possible. We have $c-b=\min\{c,a\}$ and $0\le b\le c<\mu$. Also $b\le h-a$: this is immediate for $a\ge c$, while for $a<c$ it follows from $b=c-a<h-a$. Finally,
		$$
		|D|-(\mu-b)=(q-2)\mu+\max\{0,a-c\}\ge0.
		$$
		The chosen set satisfies
		\begin{equation}\label{eq:chosen-zeros}
			|E|=(c-b)+b+b+(\mu-b)=h,\qquad
			|E\cap Z(P)|=|E\cap Z(Q)|=c.
		\end{equation}
		
		The construction is summarized in the following diagram, where shaded cells indicate the zero coordinates of the corresponding vectors.
		\begin{center}
			\begin{minipage}{\textwidth}
				\centering
				\begin{tikzpicture}[font=\small]
					\def\w{2.35}
					
					\foreach \j/\lab/\sz in {
						0/$A$/$a$,
						1/$B$/$h-a$,
						2/$C$/$h-a$,
						3/$D$/$(q-1)\mu-c+a$}
					{
						\pgfmathsetmacro\xx{\j*\w}
						\node at (\xx+0.5*\w,0.45) {\lab};
						\node[font=\scriptsize] at (\xx+0.5*\w,0.08) {size \sz};
					}
					
					\node[anchor=east] at (-0.15,-0.55) {$x$};
					\node[anchor=east] at (-0.15,-1.35) {$y$};
					\node[anchor=east] at (-0.15,-2.25) {$z$};
					
					\foreach \j/\p/\q in {
						0/0/0,
						1/0/$\ne0$,
						2/$\ne0$/0,
						3/$\ne0$/$\ne0$}
					{
						\pgfmathsetmacro\xx{\j*\w}
						
						\ifnum\j<2
						\fill[black!12] (\xx,-0.9) rectangle +(\w,0.7);
						\fi
						\ifnum\j=0
						\fill[black!12] (\xx,-1.7) rectangle +(\w,0.7);
						\fi
						\ifnum\j=2
						\fill[black!12] (\xx,-1.7) rectangle +(\w,0.7);
						\fi
						
						\draw (\xx,-0.9) rectangle +(\w,0.7);
						\draw (\xx,-1.7) rectangle +(\w,0.7);
						
						\node at (\xx+0.5*\w,-0.55) {\p};
						\node at (\xx+0.5*\w,-1.35) {\q};
					}
					
					\foreach \j/\leftlab/\rightlab in {
						0/{$E\cap A$\\[-2pt]\scriptsize $c-b$}/{$A\setminus E$},
						1/{$E\cap B$\\[-2pt]\scriptsize $b$}/{$B'$},
						2/{$E\cap C$\\[-2pt]\scriptsize $b$}/{$C'$},
						3/{$E\cap D$\\[-2pt]\scriptsize $\mu-b$}/{$D'$}}
					{
						\pgfmathsetmacro\xx{\j*\w}
						
						\fill[black!12] (\xx,-2.65) rectangle +(0.5*\w,0.8);
						\draw (\xx,-2.65) rectangle +(0.5*\w,0.8);
						\draw (\xx+0.5*\w,-2.65) rectangle +(0.5*\w,0.8);
						
						\node[align=center,font=\scriptsize]
						at (\xx+0.25*\w,-2.25) {\leftlab};
						\node[font=\scriptsize]
						at (\xx+0.75*\w,-2.25) {\rightlab};
					}
				\end{tikzpicture}
				
				\vspace{0.4em}
				
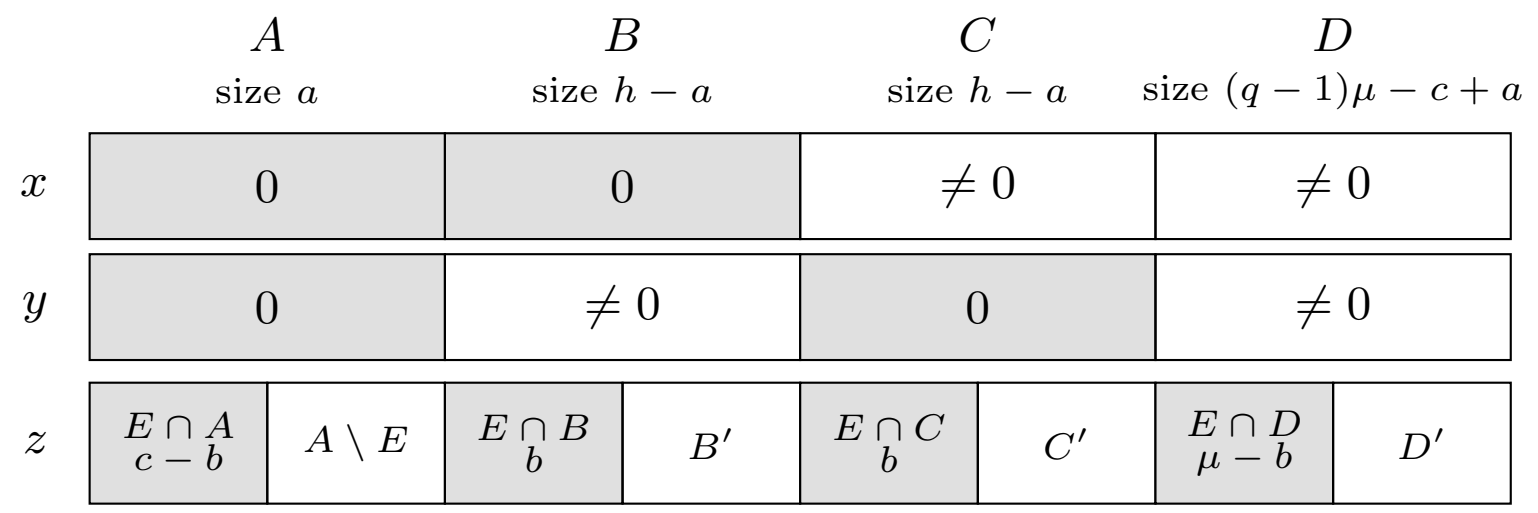
\captionof{figure}{Construction of a common neighbour for points with distinct zero sets.}
				\label{fig:distinct-zero-bridge}
			\end{minipage}
		\end{center}
		
		Write $B'=B\setminus E$, $C'=C\setminus E$ and $D'=D\setminus E$. Then
		$$
		s:=|B'|=|C'|=h-a-b=\mu-\max\{0,a-c\}.
		$$
		Since $a<h=c+\mu$, we have $1\le s\le\mu$. By \eqref{eq:blocks},
		$$
		|D'|=(q-1)\mu-s<(q-1)\mu.
		$$
		The indices in $C'$ will complete the ratios for $P$ and those in $B'$ will complete the ratios for $Q$.
		
		Suppose first that $D'\ne\varnothing$. For $i\in D'$, define
		$$
		\lambda_i=\frac{y_i}{x_i}.
		$$
		Since
		$$
		1\le |D'|<(q-1)\mu,
		$$
		Proposition~\ref{prop:Hall}(1), applied in $G=\F_q^\times$ with $r=q-1$ and $m=\mu$, gives non-zero elements $a_i$, $i\in D'$, such that every element of $\F_q^\times$ occurs at most $\mu$ times both among the $a_i$ and among the $a_i\lambda_i^{-1}$. Define
		$$
		z_i=x_i a_i.
		$$
		Then
		$$
		\frac{z_i}{x_i}=a_i,\qquad
		\frac{z_i}{y_i}=a_i\lambda_i^{-1}.
		$$
		For each $\alpha\in\F_q^\times$, let $u_\alpha$ and $v_\alpha$ be its multiplicities among the $a_i$ and among the $a_i\lambda_i^{-1}$, respectively. Then $0\le u_\alpha,v_\alpha\le\mu$ and
		$$
		\sum_{\alpha\in\F_q^\times}(\mu-u_\alpha)
		=
		\sum_{\alpha\in\F_q^\times}(\mu-v_\alpha)
		=
		(q-1)\mu-|D'|
		=
		s.
		$$
		Assign to the indices $i\in C'$ non-zero values of $z_i/x_i$ so that each $\alpha\in\F_q^\times$ occurs $\mu-u_\alpha$ times. This is possible because $|C'|=s$ and $x_i\ne0$ on $C'$. Similarly, assign to the indices $i\in B'$ non-zero values of $z_i/y_i$ so that each $\alpha\in\F_q^\times$ occurs $\mu-v_\alpha$ times.
		
		Set $z_i=0$ on $E$ and choose arbitrary non-zero $z_i$ on $A\setminus E$. Every coordinate is now defined and $Z(\Span{z})=E$, so $R=\Span{z}$ is a simplex point. By \eqref{eq:chosen-zeros},
		$$
		|Z(P)\cap Z(R)|=|Z(Q)\cap Z(R)|=c.
		$$
		On
		$$
		\supp(x)\cap\supp(z)=C'\cup D'
		$$
		each non-zero ratio occurs precisely $\mu$ times: it occurs $u_\alpha$ times on $D'$ and $\mu-u_\alpha$ times on $C'$. 
		Similarly, on
		$$
		\supp(y)\cap\supp(z)=B'\cup D',
		$$
		each non-zero ratio occurs precisely $\mu$ times. Therefore, Lemma~\ref{lem:ratio} gives
		$
		P\sim R\sim Q.
		$
		This proves the assertion when $D'\ne\varnothing$.
		
		Suppose now that $D'=\varnothing$. Then
		$$
		s=(q-1)\mu.
		$$
		Since $s\le\mu$ and $\mu>0$, it follows that $q-1\le1$. Hence $q=2$ and $s=\mu$. 
		Thus $\F_q^\times=\{1\}$. Define
		$$
		z_i=x_i\qquad\text{ for }i\in C'$$
		and
		$$
		z_i=y_i\qquad\text{ for }i\in B'.
		$$
		We set $z_i=0$ on $E$ and choose arbitrary non-zero $z_i$ on $A\setminus E$. Again $Z(\Span{z})=E$, so $R=\Span{z}$ is a simplex point. By \eqref{eq:chosen-zeros},
		$$
		|Z(P)\cap Z(R)|=|Z(Q)\cap Z(R)|=c.
		$$
		Moreover,
		$$
		\supp(x)\cap\supp(z)=C',\qquad
		\supp(y)\cap\supp(z)=B',
		$$
		and $|B'|=|C'|=s=\mu$. The only element of $\F_q^\times$ is $1$, so
		$$
		\frac{z_i}{x_i}=1\qquad\text{ for }i\in C'$$
		and
		$$
		\frac{z_i}{y_i}=1\qquad\text{ for }i\in B'.
		$$
		Hence Lemma~\ref{lem:ratio} again gives
		$
		P\sim R\sim Q.
		$
	\end{proof}
	
	\section{Points with the same zero set}\label{sec:same}
	\begin{lem}\label{lem:same-zero}
		Let $P=\Span{x}$ and $Q=\Span{y}$ be distinct simplex points with $Z(P)=Z(Q)=Z$. Define $D=\{1,\ldots,n\}\setminus Z$ and $\lambda_i=y_i/x_i$ for $i\in D$. Then $P,Q$ have a common neighbour if and only if \eqref{eq:product-condition} holds for some $F\subset D$.
	\end{lem}
	\begin{proof}
		Suppose first that $R=\Span{z}$ is collinear to both points. By Lemma~\ref{lem:ratio}, $Z(R)$ meets $Z$ in $c$ indices. Since $|Z(R)|=h=c+\mu$, the set $W=Z(R)\cap D$ has size $\mu$. Its complement $F=D\setminus W$ has size $(q-1)\mu$ and 
		$$F=\supp(x)\cap\supp(z)=\supp(y)\cap\supp(z).$$
		
		Set $a_i=z_i/x_i$ for $i\in F$. By \eqref{eq:adjacency}, each element of $\F_q^\times$ occurs precisely $\mu$ times both among the $a_i$ and among the $a_i\lambda_i^{-1}=z_i/y_i$. Hence the two sequences have equal products:
		$$
		\prod_{i\in F}a_i
		=\left(\prod_{\alpha\in\F_q^\times}\alpha\right)^\mu
		=\prod_{i\in F}a_i\lambda_i^{-1}.
		$$
		Cancelling $\prod_{i\in F}a_i$ gives $\prod_{i\in F}\lambda_i=1$.
		
		Conversely, let $F\subset D$ satisfy \eqref{eq:product-condition} and define $W=D\setminus F$. 
		Proposition~\ref{prop:Hall}(2), applied in $G=\F_q^\times$ with
		$r=q-1$ and $m=\mu$, gives elements $a_i\in\F_q^\times$, $i\in F$, such that every non-zero field element occurs precisely $\mu$ times both among the $a_i$ and among the $a_i\lambda_i^{-1}$. Choose $Z_0\subset Z$ of size $c$ and define a vector $z$ by
		$$
		z_i=
		\begin{cases}
			0,&i\in W\cup Z_0,\\
			x_i a_i,&i\in F,\\
			\text{an arbitrary element of }\F_q^\times,&i\in Z\setminus Z_0.
		\end{cases}
		$$
		The construction is summarized in the following diagram, where shaded cells indicate zero coordinates.
		\begin{center}
			\begin{minipage}{\textwidth}
				\centering
				\begin{tikzpicture}[font=\small]
					\def\w{2.3}
					
					\foreach \j/\lab/\sz in {
						0/$Z_0$/$c$,
						1/$Z\setminus Z_0$/$\mu$,
						2/$W$/$\mu$,
						3/$F$/$(q-1)\mu$}
					{
						\pgfmathsetmacro\xx{\j*\w}
						\node at (\xx+0.5*\w,0.45) {\lab};
						\node[font=\scriptsize] at (\xx+0.5*\w,0.08) {size \sz};
					}
					
					\node[anchor=east] at (-0.15,-0.55) {$x$};
					\node[anchor=east] at (-0.15,-1.35) {$y$};
					\node[anchor=east] at (-0.15,-2.15) {$z$};
					
					\foreach \j/\val in {
						0/0,
						1/0,
						2/$\ne0$,
						3/$\ne0$}
					{
						\pgfmathsetmacro\xx{\j*\w}
						
						\ifnum\j<2
						\fill[black!12] (\xx,-0.9) rectangle +(\w,0.7);
						\fi
						
						\draw (\xx,-0.9) rectangle +(\w,0.7);
						\node at (\xx+0.5*\w,-0.55) {\val};
					}
					
					\foreach \j/\val in {
						0/0,
						1/0,
						2/$\ne0$,
						3/$\ne0$}
					{
						\pgfmathsetmacro\xx{\j*\w}
						
						\ifnum\j<2
						\fill[black!12] (\xx,-1.7) rectangle +(\w,0.7);
						\fi
						
						\draw (\xx,-1.7) rectangle +(\w,0.7);
						\node at (\xx+0.5*\w,-1.35) {\val};
					}
					
					\foreach \j/\val in {
						0/0,
						1/$\ne0$,
						2/0,
						3/$\ne0$}
					{
						\pgfmathsetmacro\xx{\j*\w}
						
						\ifnum\j=0
						\fill[black!12] (\xx,-2.5) rectangle +(\w,0.7);
						\fi
						\ifnum\j=2
						\fill[black!12] (\xx,-2.5) rectangle +(\w,0.7);
						\fi
						
						\draw (\xx,-2.5) rectangle +(\w,0.7);
						\node at (\xx+0.5*\w,-2.15) {\val};
					}
				\end{tikzpicture}
				
				\vspace{0.4em}
				\captionof{figure}{Construction of a common neighbour for points with the same zero set.}
				\label{fig:same-zero-bridge}
			\end{minipage}
		\end{center}
		The zero set of $R=\Span{z}$ has size $|W|+|Z_0|=\mu+c=h$, so $R$ is a simplex point. 
		For both pairs $(P,R)$ and $(Q,R)$, the intersection of the zero sets is
		$Z_0$ and the intersection of the supports is $F$.
		The two sequences given by Proposition~\ref{prop:Hall}(2) have the required ratio multiplicities. 
		The point $R$ is distinct from both $P,Q$ since $|Z\cap Z(R)|=c<h$ and Lemma~\ref{lem:ratio} gives $P\sim R\sim Q$.
	\end{proof}
	
	The equivalent condition using $W$ follows by taking complements:
	$$
	\prod_{i\in F}\lambda_i=1
	\quad\text{if and only if}\quad
	\prod_{i\in W}\lambda_i=\prod_{i\in D}\lambda_i=\Lambda.
	$$
	
	\begin{rem}\label{rem:k2}{\rm
			For $k=2$ we have $c=0$ and $\mu=1$. 
			The set $W$ consists of one index and the condition becomes
			\begin{equation}\label{eq:k2-condition}
				P,Q\text{ have a common neighbour if and only if }
				\Lambda\in\{\lambda_i:i\in D\}.
			\end{equation}
			More generally, for points with the same zero set we have
			$|D|=q^{k-1}\equiv1\pmod{{q-1}}$. 
			Hence replacing $y$ by
			$\Lambda^{-1}y$ makes the product of all ratios equal to one.
			For $k=2$, in this normalization \eqref{eq:k2-condition} says that one
			of the ratios is equal to one.
	}\end{rem}
	
	\section{Completion of the proofs}\label{sec:proofs}
	\subsection{The product condition in higher dimension}
	The multiplicative group $\F_q^\times$ is cyclic. 
	We use the classical Erd\H{o}s--Ginzburg--Ziv theorem \cite{EGZ61}, stated here in multiplicative form to apply directly to the coordinate ratios.
	
	\begin{theorem}[Erd\H{o}s--Ginzburg--Ziv]\label{thm:EGZ}
		Let $G$ be a cyclic group of order $r$. Every sequence of $2r-1$ elements of $G$ contains $r$ terms with distinct indices whose product is the identity.
	\end{theorem}
	
	\begin{lem}\label{lem:iterated-EGZ}
		Let $G$ be a cyclic group of order $r$ and let $t\ge1$ be an integer. 
		Every sequence of at least $rt+r-1$ elements of $G$ contains $rt$ terms with distinct indices whose product is the identity.
	\end{lem}
	\begin{proof}
		Choose and remove $r$ terms with product one, repeating this $t$ times. 
		Before the $j$-th choice, where $1\le j\le t$, there are at least
		$$
		rt+r-1-(j-1)r=r(t-j+1)+r-1\ge2r-1
		$$
		terms left. Therefore, Theorem~\ref{thm:EGZ} applies at every step. The union of the $t$ disjoint selections has size $rt$ and product one. The case $r=1$ is immediate and is also covered by this argument.
	\end{proof}
	
	\begin{prop}\label{prop:automatic}
		If $\mu=q^{k-2}\ge q-2$, then any two distinct simplex points with the same zero set have a common neighbour. In particular, this holds for $k\ge3$ and for $(k,q)=(2,3)$.
	\end{prop}
	\begin{proof}
		The sequence of ratios in Lemma~\ref{lem:same-zero} has length $q\mu$. Apply Lemma~\ref{lem:iterated-EGZ} to $\F_q^\times$ with $r=q-1$ and $t=\mu$. Since
		$$
		rt+r-1=(q-1)\mu+q-2\le q\mu,
		$$
		we obtain $(q-1)\mu$ ratios whose product is one. Therefore, Lemma~\ref{lem:same-zero} gives a common neighbour. If $k\ge3$, then $\mu=q^{k-2}\ge q\ge q-2$. For $(k,q)=(2,3)$, the required inequality is $1\ge1$.
	\end{proof}
	
	\subsection{Distances and diameters}
	\begin{proof}[Proof of Theorem~\ref{thm:distance}]
		Part (1) is Lemma~\ref{lem:ratio}. In part (2), Lemma~\ref{lem:bridge} gives a common neighbour, so the distance is two because $P,Q$ are distinct and non-adjacent. If the zero sets coincide, the points are not adjacent. Therefore, Lemma~\ref{lem:same-zero} gives precisely the distance-two assertion in part (3) and the equivalent formulation in terms of $W$ is proved after that lemma.
		
		It remains to determine the distance when the product condition fails. Thus, the points have no common neighbour and no path of length at most two joins them. Choose a neighbour $R$ of $P$. 
		Lemma~\ref{lem:ratio} gives
		$$
		|Z(P)\cap Z(R)|=c<h,
		$$
		hence $Z(R)\ne Z(P)=Z(Q)$. Lemma~\ref{lem:bridge}, applied to $R,Q$, gives a point $S$ with $R\sim S\sim Q$. Thus $P\sim R\sim S\sim Q$ is a walk of length three and $d_\Delta(P,Q)=3$. Every pair is joined by a path, which proves connectedness.
	\end{proof}
	
	\begin{proof}[Proof of Theorem~\ref{thm:diameter}]
		If $k\ge3$ or $(k,q)=(2,3)$, Lemma~\ref{lem:bridge} and Proposition~\ref{prop:automatic} give a common neighbour for any two distinct vertices. Therefore, the diameter is at most two. 
		To show that the graph is not complete, we first consider the case $q>2$. 
		Fix a zero set of size $h$ and choose two vectors with that zero set which differ in precisely one non-zero coordinate. Their support has size $q^{k-1}\ge2$, so the vectors are not proportional. The corresponding simplex points are distinct and not collinear.
		
		For $q=2$ and $k\ge3$, take the zero sets $\{1,\ldots,h\}$ and $\{1,\ldots,h-1,h+1\}$. They determine simplex points and their intersection has size $h-1$ with
		$$
		h-1-c=\mu-1>0.
		$$
		Since $h-1>c$, the first condition in \eqref{eq:adjacency} fails, so these points are not adjacent. 
		Thus the diameter is two in all the stated cases.
		
		If $k=2$ and $q=2$, the simplex points are represented by $(0,1,1)$, $(1,0,1)$ and $(1,1,0)$. 
		Every pair satisfies \eqref{eq:adjacency}, so $\Delta(2,2)$ is the complete graph on three vertices and its diameter is one.
		
		Finally, let $k=2$ and $q\ge4$. Theorem~\ref{thm:distance} gives the upper bound three. Choose a generator $a$ of the cyclic group $\F_q^\times$ and consider the points represented by
		$$
		x=(0,1,\ldots,1),\qquad
		y=(0,a,a,1,\ldots,1)
		$$
		in $\F_q^{q+1}$. Both have weight $q$ and the same zero set. Their $q$ ratios consist of two copies of $a$ and $q-2$ copies of $1$, so $\Lambda=a^2$. Since the order of $a$ is $q-1\ge3$, we have $a^2\notin\{1,a\}$. Condition~\eqref{eq:k2-condition} fails and the two points are at distance three. This proves the lower bound and completes the diameter calculation.
	\end{proof}
	
	\subsection*{Acknowledgment}
	The author is grateful to Mark Pankov for useful discussions.

\end{document}